\documentclass[preprint,12pt]{elsarticle}

\usepackage{amssymb}
\usepackage{amsmath}
\usepackage{amsthm}
\usepackage{enumerate}

\newtheorem{defn}{Definition}
\newtheorem{prop}{Proposition}
\newtheorem{exmp}{Example}
\newtheorem{rem}{Remark}

\newproof{pf}{Proof}
\newproof{ack}{Acknowledgements}

\journal{}

\begin{document}

\begin{frontmatter}


\cortext[cor1]{Corresponding author (Fax:+903562521585)}

\title{Neutrosophic soft sets and neutrosophic soft matrices based on decision making}


\author[rvt1]
{Irfan Deli\corref{cor1}} \ead{irfandeli@kilis.edu.tr}
\author[rvt3]{{Said Broumi}}
\ead{broumisaid78@gmail.com }
\address[rvt1]{ 7 Aral{\i}k
University, 79000 Kilis, Turkey}
\address [rvt3]{ Administrator of Faculty of Arts and Humanities, Hay El Baraka Ben M'sik Casablanca B.P. 7951, Hassan II
University Mohammedia-Casablanca , Morocco }
\begin{abstract}

Maji\cite{maj-13}, firstly proposed neutrosophic soft sets can
handle the indeterminate information and inconsistent information
which exists commonly in belief systems.  In this paper, we have
firstly redefined complement, union and compared our definitions of
neutrosophic soft with the definitions given by Maji. Then, we have
introduced the concept of neutrosophic soft matrix and their
operators which are more functional to make theoretical studies in
the neutrosophic soft set theory. The matrix is useful for storing
an neutrosophic soft set in computer memory which are very useful
and applicable. Finally, based on some of these matrix operations a
efficient methodology named as NSM-decision making has been
developed to solve neutrosophic soft set based group decision making
problems.
\end{abstract}

\end{frontmatter}


\section{Introduction}
In recent years a number of theories have been proposed to deal with
uncertainty, imprecision, vagueness and indeterminacy. Theory of
probability, fuzzy set theory\cite{zad-65}, intuitionistic fuzzy
sets \cite{ata-99}, interval valued intuitionistic fuzzy sets
\cite{ata-89}, vague sets\cite{gau-93}, rough set
theory\cite{paw-82}, neutrosophic theory\cite{sam-05}, interval
neutrosophic theory\cite{wan-05} etc. are consistently being
utilized as efficient tools for dealing with diverse types of
uncertainties and imprecision embedded in a system. However, each of
these theories has its inherent difficulties as pointed out by
Molodtsov\cite{mol-99}. The reason for these difficulties is,
possibly, the inadequacy of the parameterization tool of the
theories. Later on, many interesting results of soft set theory have
been obtained by embedding the idea of fuzzy set, intuituionstic
fuzzy set, vague set, rough set, interval intuitionistic fuzzy set,
intuitionistic neutrosophic set, interval neutrosophic set,
neutrosophic set and so on. For example, fuzzy soft
set\cite{maj-01a}, intuitionistic fuzzy soft
set\cite{cag-13,maj-01b}, rough soft set\cite{fen-11,fen-10},
interval valued intuitionistic fuzzy soft
set\cite{jia-10,yan-09,zha-13}, neutrosophic soft set\cite{maj-13,
maj-12}, generalized neutrosophic soft set\cite{bro-13},
intuitionstic neutrosophic soft set\cite{bro-13b}, interval valued
neutrosophic soft set\cite{del-14}. The theories has developed in
many directions and applied to wide variety of fields such as  on
soft decision making\cite{cag-09ss,zou-08}, fuzzy soft decision
making\cite{deli-12a,deli-12b,kon-09,roy-07} ,on relation of fuzzy
soft set\cite{som-06,süt-12}, on relation on intuiotionstic fuzzy
soft set\cite{din-10,muk-08}, on relation on neutrosophic soft
set\cite{del-14c}, on relation on interval neutrosophic soft
set\cite{del-14} and so on.

Researchers published several papers on fuzzy soft matrices and
intuitionistic fuzzy soft matrices, and it has been applying in many
fields of real life
scenarios(see\cite{kal-13,kon-09,kha-13,mao-13}). Recently Cagman et
al\cite{cag-09sm} introduced soft matrices and applied it in
decision making problem. They also introduced introduced fuzzy soft
matrices\cite{cag-12}, Chetia and Das\cite{che-x} defined
intuitionistic fuzzy soft matrices with different products and
properties on these products. Further Saikia et al\cite{sai-14}
defined generalized fuzzy soft matrices with four different product
of generalized intuitionstic fuzzy soft matrices and presented an
application in medical diagnosis. Next, Broumi et al\cite{bro-13c}
studied fuzzy soft matrix based on reference function and defined
some new operations such fuzzy soft complement matrix, trace of
fuzzy soft matrix based on reference function a new fuzzy soft
matrix decision method based on reference function is presented.
Recently, Mondal et al\cite{mon-13a,mon-14,mon-13b} introduced fuzzy
and intuitionstic fuzzy soft matrix and multicrita in desion making
based on three basic t-norm operators. The matrices has differently
developed in many directions and applied to wide variety of fields
in \cite{bas-13,bor-12,raj-13,sai-13}.

Our objective is to introduce the concept of neutrosophic matrices
and its applications in decision making problem. The remaining part
of this paper is organized as follows. Section 2 contains basic
definitions and notations that are used in the remaining parts of
the paper. we investigated redefined neutrosophic soft set and some
operations and compared our definitions of neutrosophic soft with
the definitions given
 Maji\cite{maj-13} in section 3. In section 4, we introduce the concept of
neutrosophic matrices and presente some of theirs basic properties.
In section 5, we present two special products of neutrosophic
matrices. In section 6, we present a soft decision making method
based on and-product of neutrsophic matrice. Finally ,conclusion is
made in section 7.

\section{Preliminary}\label{ss}
In this section, we give the basic definitions and results of
neutrosophic set theory \cite{sam-05}, soft set theory \cite{mol-99}
and soft matrix theory \cite{cag-09sm} that are useful for
subsequent discussions.

\begin{defn} \cite{sam-05}
Let U be a space of points (objects), with a generic element in U
denoted by u. A neutrosophic sets(N-sets) A in U is characterized by
a truth-membership function $T_A$, a indeterminacy-membership
function $I_A$ and a falsity-membership function $F_A$. $T_A(u)$;
$I_A(u)$ and $F_A(u)$ are real standard or nonstandard subsets of
$[0,1]$. It can be written as

$$A=\{<u,(T_A(u),I_A(u),F_A(u))>:u\in U, \,T_A(u),I_A(u),F_A(u)\in [0,1]\}.$$

There is no restriction on the sum of $T_A(u)$; $I_A(u)$ and
$F_A(u)$, so $0\leq sup T_A(u) + sup I_A(u) + supF_A(u)\leq 3$.
\end{defn}
\begin{defn}\cite{mol-99}
Let $U$ be a universe, $E$ be a set of parameters that are describe
the elements of $U$, and $A\subseteq E$. Then, a soft set $F_A$ over
$U$ is a set defined by a set valued function $f_A$ representing a
mapping
\begin{equation}\label{soft-set}
f_A: E\to P(U) \textrm{ such that}\, f_A(x)=\emptyset \textrm{ if }
x\in E-A
\end{equation}
where $f_A$ is called approximate function of the soft set $F_A$. In
other words, the soft set is a parametrized family of subsets of the
set $U$, and therefore it can be written a set of ordered pairs
$$
F_A= \{(x, f_A(x)): x\in E, f_A(x)=\emptyset \textrm{ if } x\in
E-A\}
$$
\end{defn}
The subscript $A$ in the $f_A$ indicates that $f_A$ is the
approximate function of $F_A$. The value $f_A(x)$ is a set called
\emph{$x$-element} of the soft set for every $x\in E$.
\begin{defn}\label{matrix}\cite{cag-09sm}
Let $F_A$ be a soft set over U. Then a subset of $U\times E$ is
uniquely defined by
$$
R_A=\{((u,x) /(u,x)): (u,x)\in U\times E\}
$$
which is called a relation form of $F_A$. The characteristic
function of $R_A$ is written by
\begin{displaymath}
\chi_{R_A}:U\times E\to[0,1],\quad \chi_{R_A}(u,x) = \left\{
  \begin{array}{ll}
   1,&\hbox{if } (u,x)\in R_A,\\
    0,&\hbox{if }  (u,x)\notin R_A.
  \end{array}
\right.
\end{displaymath}

If $U=\{u_1,u_2,...,u_m\}$, $E=\{x_1,x_2,...,x_n\}$ and $A\subseteq
E$, then the $R_A$ can be presented by a table as in the following
form
$$
\begin{tabular}{c|cccc}
\({R_A}\)& \(x_{1}\)                  & \(x_{2}\)                  &...&\(x_{n}\) \\
\hline
\(u_1\)  & \(\chi_{R_A}(u_{1},x_{1})\)& \(\chi_{R_A}(u_{1},x_{2})\)&...& \(\chi_{R_A}(u_{1},x_{n})\) \\
\(u_2\)  & \(\chi_{R_A}(u_{2},x_{1})\)& \(\chi_{R_A}(u_{2},x_{2})\)&...& \(\chi_{R_A}(u_{2},x_{n})\) \\
\vdots   & \vdots                     &  \vdots                &\(\ddots\)&  \vdots \\
\(u_m\)  & \(\chi_{R_A}(u_{m},x_{1})\)& \(\chi_{R_A}(u_{m},x_{2})\)&...& \(\chi_{R_A}(u_{m},x_{n})\) \\
\end{tabular}
$$
If $a_{ij}=\chi_{R_A}(u_i,x_j)$, we can define a matrix
$$
[a_{ij}]_{m\times n}=\left[
\begin{array}{cccc}
x_{11} & x_{12} & \cdots & x_{1n}\\
x_{21} & x_{22} & \cdots & x_{2n}\\
\vdots & \vdots &  \ddots &    \vdots\\
x_{m1} & x_{m2}& \cdots & x_{mn}\\
\end{array}\right ]
$$
which is called an $m\times n$ $s$-matrix of  the soft set $F_A$
over $U$.
\end{defn}
From now on we shall delete the subscripts $m\times n$ of
$[a_{ij}]_{m\times n}$, we use $[a_{ij}]$ instead of
$[a_{ij}]_{m\times n}$ for $i=1,2,...m$ and $j=1,2,...n$.
\begin{defn}\label{and}\cite{cag-09sm}
Let $[a_{ij}],[b_{ik}]\in FSM_{m\times n}$. Then $And$-product of
$[a_{ij}]$ and $[b_{ik}]$ is defined by
$$
\wedge:FSM_{m\times n}\times FSM_{m\times n}\to FSM_{m\times
n^2},\quad [a_{ij}]\wedge[b_{ik}]=[c_{ip}]
$$
where  $c_{ip}=\min\{a_{ij}, b_{ik}\}$ such that $p=n(j-1)+k$.
\end{defn}
%

\begin{defn}\label{and}\cite{cag-09sm}
Let $[a_{ij}],[b_{ik}]\in FSM_{m\times n}$. Then $Or$-product of
$[a_{ij}]$ and $[b_{ik}]$ is defined by
$$
\vee:FSM_{m\times n}\times FSM_{m\times n}\to FSM_{m\times
n^2},\quad [a_{ij}]\vee[b_{ik}]=[c_{ip}]
$$
where  $c_{ip}=\max\{a_{ij}, b_{ik}\}$ such that $p=n(j-1)+k$.
\end{defn}
\begin{defn}\label{SMm}\cite{cag-09sm}
Let $[c_{ip}]\in SM_{m\times n^2}$, $I_k=\{p: \exists i, c_{ip}\neq
0, (k-1)n<p\leqslant kn\}$ for all $k\in I=\{1,2,...,n\}$. Then
$fs$-max-min decision function, denoted $Mm$, is defined as follows
$$
Mm:FSM_{m\times n^2}\to FSM_ {m\times 1},\quad
Mm[c_{ip}]=[d_{i1}]=[\max_{k}\{t_{ik}\}]
$$
where
$$
t_{ik}= \left\{
  \begin{array}{ll}
    \min_{p\in I_k}\{c_{ip}\},&\hbox{if }  I_k\neq \emptyset,\\
    0,&\hbox{if }  I_k=\emptyset.
  \end{array}
\right.
$$
The one column $fs$-matrix $Mm[c_{ip}]$ is called max-min decision
$fs$-matrix.
\end{defn}
\begin{defn}\label{SMm}
Let $U=\{u_1, u_2,...,u_m\}$ be an initial universe and
$Mm[c_{ip}]=[d_{i1}]$. Then a subset of $U$ can be obtained by using
$[d_{i1}]$ as in the following way
$$opt_{[d_{i1}]}(U)=\{d_{i1}/u_i:u_i\in U, d_{i1}\neq 0\}$$
which is called an optimum fuzzy set on $U$.
\end{defn}
\begin{defn}\cite{maj-13}
Let $U$ be a universe, $N(U)$ be the set of all neutrosophic sets on
U, $E$ be a set of parameters that are describe the elements of $U$,
and $A\subseteq E$. Then, a neutrosophic soft set $N$ over $U$ is a
set defined by a set valued function $f_N$ representing a mapping
$$f_N: A\to N(U) $$
where $f_N$ is called approximate function of the neutrosophic soft
set $N$. In other words, the neutrosophic soft set is a parametrized
family of some elements of the set $P(U)$, and therefore it can be
written a set of ordered pairs
$$
N= \{(x, f_N(x)): x\in A\}
$$
\end{defn}

\begin{defn}\cite{maj-13}
Let $N_1$ and $N_2$ be two neutrosophic soft sets over neutrosophic
soft universes $(U,A)$ and $(U, B)$, respectively.
\begin{enumerate}
    \item $N_1$ is said to be neutrosophic soft subset of $N_2$ if $A\subseteq
B$ and $T_{f_{N_1(x)}}(u)\leq T_{f_{N_2(x)}}(u)$,
$I_{f_{N_1(x)}}(u)\leq I_{f_{N_2(x)}}(u)$ ,$F_{f_{N_1(x)}}(u)\geq
F_{f_{N_2(x)}}(u)$,  $\forall x\in A$, $u\in U$.
    \item $N_1$ and $N_2$ are said to be equal if $N_1$ neutrosophic soft subset of $N_2$
   and $N_2$ neutrosophic soft subset of $N_2$.
\end{enumerate}
\end{defn}

\begin{defn}\cite{maj-13}
Let $E = \{e_1, e_2,...\}$ be a set of parameters. The NOT set of E
is denoted by $\neg E$ is defined by $\neg E=\{\neg e_1, \neg
e_2,...\}$ where $\neg e_i= not \, e_i,\forall i$.
\end{defn}
\begin{defn}\cite{maj-13}
Let $N_1$ and $N_2$ be two neutrosophic soft sets over soft
universes $(U,A)$ and $(U, B)$, respectively,
\begin{enumerate}
    \item The complement of a neutrosophic soft set $N_1$ denoted by $N_1^{\circ}$
and is defined by a set valued function $f_N^{\circ}$ representing a
mapping $f_{N_1}^{\circ}: \neg E\to N(U)$

$f_{N_1}^{\circ}=\{(u,
<F_{f_{N_1(x)}}(u),I_{f_{N_1(x)}}(u),T_{f_{N_1(x)}}(u)>): x\in \neg
E, u\in U\}.$

    \item Then the union of $N_1$ and $N_2$ is denoted by $N_1\acute{{\cup}} N_2$ and
    is defined by $N_3$$(C=A\cup B)$, where the truth-membership, indeterminacy-membership and
falsity-membership of $N_3$ are as follows: $\forall u\in
  U$,

$$
T_{f_{N_3(x)}}(u)=\left \{\begin{array}{ll}
  \displaystyle T_{f_{N_1(x)}}(u) , &if x\in A-B
  \\
  \displaystyle T_{f_{N_2(x)}}(u) , & if x\in B- A \\
  \displaystyle max\{T_{f_{N_1(x)}}(u),T_{f_{N_2(x)}}(u)\} , & if x\in A\cap B
\end{array}\right.
$$

$$
I_{f_{N_3(x)}}(u)=\left \{\begin{array}{ll}
  \displaystyle I_{f_{N_1(x)}} (u), &if x\in A-B
  \\
  \displaystyle I_{f_{N_2(x)}}(u) , & if x\in B- A \\
  \displaystyle \frac{(I_{f_{N_1(x)}}(u),I_{f_{N_2(x)}}(u))}{2} , & if x\in A\cap B
\end{array}\right.
$$

$$
F_{f_{N_3(x)}}(u)=\left \{\begin{array}{ll}
  \displaystyle F_{f_{N_1(x)}}(u) , &if x\in A-B
  \\
  \displaystyle F_{f_{N_2(x)}}(u) , & if x\in B- A \\
  \displaystyle min\{I_{f_{N_1(x)}}(u),I_{f_{N_2(x)}}(u)\} , & if x\in A\cap B
\end{array}\right.
$$

    \item Then the intersection of $N_1$ and $N_2$ is denoted by $N_1 \acute{{\cap}} N_2$ and
    is defined by $N_3$$(C=A\cap B)$, where the truth-membership, indeterminacy-membership and
falsity-membership of $N_3$ are as follows: $\forall u\in
  U$,

$
 T_{f_{N_3(x)}}(u)=min\{T_{f_{N_1(x)}}(u),T_{f_{N_2(x)}}(u)\}$
  , $I_{f_{N_3(x)}}(u)=\frac{(I_{f_{N_1(x)}}(u),I_{f_{N_2(x)}}(u))}{2}$

  and $F_{f_{N_3(x)}}(u)=max\{F_{f_{N_1(x)}}(u),F_{f_{N_2(x)}}(u)\},\, \forall x\in
  C.
$
\begin{defn}\label{add}\cite{dub-80}
$t$-norms are associative, monotonic and commutative two valued
functions $t$ that map from $[0,1]\times [0,1]$ into $[0,1]$. These
properties are formulated with the following conditions:  $\forall
a,b,c,d \in [0,1],$
\begin{enumerate}[{i.}]
\item $t(0,0)=0$ and $t(a,1)=t(1,a)=a,$
\item If  $a\leq c$ {and} $b\leq d$, then $t(a,b)\leq t(c,d)$
\item $t(a,b)= t(b,a)$
\item $t(a,t(b,c))=t(t(a,b,c))$
\end{enumerate}
\end{defn}
\begin{defn}\label{add}\cite{dub-80}
$t$-conorms ($s$-norm) are associative, monotonic and commutative
two placed functions $s$ which map from $[0,1]\times [0,1]$ into
$[0,1]$. These properties are formulated with the following
conditions: $\forall a,b,c,d \in [0,1],$
\begin{enumerate}[{i.}]
\item  $s(1,1)=1 \,\,and \,\, s(a,0)=s(0,a)=a,$
\item  {if}
$a\leq c$ {and}  $b\leq d$, then $s(a,b)\leq s(c,d)$
\item   $s(a,b)= s(b,a)$
\item  $s(a,s(b,c))=s(s(a,b,c)$
\end{enumerate}
\end{defn}
$t$-norm and $t$-conorm are related in a sense of lojical duality.
Typical dual pairs of non parametrized $t$-norm and $t$-conorm are
complied below:

\begin{enumerate}[{i.}]
\item  Drastic product:
$$
t_w(a,b) = \left \{\begin{array}{ll}
min\{a,b\}, & max\{a  b\}=1\\
0,                                & otherwise \\
\end{array}\right.
$$
\item  Drastic sum:
$$
s_w(a,b) = \left \{\begin{array}{ll}
max\{a,b\}, &  min\{a  b\}=0\\
1, & otherwise\\
\end{array}\right.
$$
\item  Bounded product:
$$
t_1(a,b)=max\{0, a+b-1\}
$$
\item  Bounded sum:
$$
s_1(a,b)=min\{1, a+b\}
$$
\item Einstein product:
 $$
 t_{1.5}(a,b)=\frac{a.b}{2-[a+b-a.b]}
 $$
\item  Einstein sum:
$$
s_{1.5}(a,b)=\frac{a+b}{1+a.b}
$$
\item  Algebraic product:
$$
t_{2}(a,b)=a.b
$$
\item  Algebraic sum:
$$
s_{2}(a,b)=a+b-a.b
$$
\item  Hamacher product: $$
t_{2.5}(a,b)=\frac{a.b}{a+b-a.b}
$$
\item  Hamacher sum:
$$
s_{2.5}(a,b)=\frac{a+b-2.a.b}{1-a.b}
$$
\item  Minumum:
$$
t_3(a,b)=min\{a,b\}
$$
\item  Maximum:
$$
s_3(a,b)=max\{a,b\}
$$
\end{enumerate}

\end{enumerate}
\end{defn}
\section{Neutrosophic soft set and some operations redefined}
Notion of the neutrosophic soft set theory is first given by Maji
\cite{maj-13}. This section, we has modified the definition of
neutrosophic soft set and operations as follows. Some of it is
quoted from \cite{ans-13,ash-02,cag-14,din-10, maj-13,sam-05}.
\begin{defn}
Let $U$ be a universe, $N(U)$ be the set of all neutrosophic sets on
U, $E$ be a set of parameters that are describe the elements of $U$
Then, a neutrosophic soft set $N$ over $U$ is a set defined by a set
valued function $f_N$ representing a mapping
$$f_N: E\to N(U)$$
where $f_N$ is called approximate function of the neutrosophic soft
set $N$. For $x\in E$, the set $f_N(x)$ is called x-approximation of
the neutrosophic soft set $N$ which may be arbitrary, some of them
may be empty and some may have a nonempty intersection. In other
words, the neutrosophic soft set is a parametrized family of some
elements of the set $N(U)$, and therefore it can be written a set of
ordered pairs,
$$
N= \{(x, \{<u,T_{f_N(x)}(u),I_{f_N(x)}(u),F_{f_N(x)}(u)>:x\in U\}:
x\in E\}
$$
where

$T_{f_N(x)}(u),I_{f_N(x)}(u),F_{f_N(x)}(u)\in [0,1]$
\end{defn}

\begin{defn}
Let $N_1$ and $N_2$ be two neutrosophic soft sets. Then, the
complement of a neutrosophic soft set $N_1$ denoted by $N_1^c$ and
is defined by
$$
 {N_1}^c= \{(x,
\{<u,F_{f_{N_1(x)}}(u),I_{f_{N_1(x)}}(u),T_{f_{N_1(x)}}(u)>:x\in
U\}: x\in E\}
$$
\end{defn}
\begin{defn}
Let $N_1$ and $N_2$ be two neutrosophic soft sets. Then, the union
of $N_1$ and $N_2$ is denoted by $N_3=N_1\tilde{\cup} N_2$ and
    is defined by
$$
 {N_3}= \{(x,
\{<u,T_{f_{N_3(x)}}(u),I_{f_{N_3(x)}}(u),F_{f_{N_3(x)}}(u)>:x\in
U\}: x\in E\}
$$
where

$ T_{f_{N_3(x)}}(u)=s(T_{f_{N_1(x)}}(u),T_{f_{N_2(x)}}(u)) $, $
I_{f_{N_3(x)}}(u)=t(I_{f_{N_1(x)}}(u),I_{f_{N_2(x)}}(u)) $ and $
F_{f_{N_3(x)}}(u)=t(F_{f_{N_1(x)}}(u),F_{f_{N_2(x)}}(u)) $
\end{defn}

\begin{defn}
Let $N_1$ and $N_2$ be two neutrosophic soft sets. Then, the
intersection of $N_1$ and $N_2$ is denoted by $N_4=N_1\tilde{\cap}
N_2$ and
    is defined by
$$
 {N_4}= \{(x,
\{<u,T_{f_{N_4(x)}}(u),I_{f_{N_4(x)}}(u),F_{f_{N_4(x)}}(u)>:x\in
U\}: x\in E\}
$$
where

$ T_{f_{N_4(x)}}(u)=t(T_{f_{N_1(x)}}(u),T_{f_{N_2(x)}}(u)) $, $
I_{f_{N_4(x)}}(u)=s(I_{f_{N_1(x)}}(u),I_{f_{N_2(x)}}(u)) $ and $
F_{f_{N_4(x)}}(u)=s(F_{f_{N_1(x)}}(u),F_{f_{N_2(x)}}(u)) $
\end{defn}

\begin{exmp}\label{1}
Let $U=\{u_1, u_2, u_3,u_4\}$, $E=\{x_1, x_2, x_3\}$. $N_1$ and
$N_2$ be two neutrosophic soft sets as
$$
\begin{array}{rr}
N_1= &
\bigg\{(x_1,\{<u_1,(0.4,0.5,0.8)>,<u_2,(0.2,0.5,0.1)>,<u_3,(0.3,0.1,0.4)>,\\&
<u_4,(0.4,0.7,0.7)>\}),
(x_2,<u_1,(0.5,0.7,0.7)>,<u_2,(0.3,0.6,0.3)>,\\&<u_3,(0.2,0.6,0.5)>,
<u_4,(0.4,0.5,0.5)>\}),
(x_3,\{<u_1,(0.7,0.8,0.6)>,\\&<u_2,(0.5,0.6,0.7)>,<u_3,(0.7,0.5,0.8)>,
<u_4,(0.2,0.8,0.5)>\})\bigg\} \\
\end{array}
$$
and
$$
\begin{array}{rr}
N_2= &
\bigg\{(x_1,\{<u_1,(0.7,0.6,0.7)>,<u_2,(0.4,0.2,0.8)>,<u_3,(0.9,0.1,0.5)>,\\&
<u_4,(0.4,0.7,0.7)>\}),
(x_2,<u_1,(0.5,0.7,0.8)>,<u_2,(0.5,0.9,0.3)>,\\&<u_3,(0.5,0.6,0.8)>,
<u_4,(0.5,0.8,0.5)>\}),
(x_3,\{<u_1,(0.8,0.6,0.9)>,\\&<u_2,(0.5,0.9,0.9)>,<u_3,(0.7,0.5,0.4)>,
<u_4,(0.3,0.5,0.6)>\})\bigg\} \\
\end{array}
$$
here;
$$
\begin{array}{rr}
N_1^c= &
\bigg\{(x_1,\{<u_1,(0.8,0.5,0.4)>,<u_2,(0.1,0.5,0.2)>,<u_3,(0.4,0.1,0.3)>,\\&
<u_4,(0.7,0.7,0.4)>\}),
(x_2,<u_1,(0.7,0.7,0.5)>,<u_2,(0.3,0.6,0.3)>,\\&<u_3,(0.5,0.6,0.2)>,
<u_4,(0.5,0.5,0.4)>\}),
(x_3,\{<u_1,(0.6,0.8,0.7)>,\\&<u_2,(0.7,0.6,0.5)>,<u_3,(0.8,0.5,0.7)>,
<u_4,(0.5,0.8,0.2)>\})\bigg\} \\
\end{array}
$$
Let us consider the t-norm $min\{a,b\}$ and s-norm $max\{a,b\}$.
Then,
$$
\begin{array}{rr}
N_1\tilde{\cup} N_2= &
\bigg\{(x_1,\{<u_1,(0.7,0.5,0.7)>,<u_2,(0.4,0.2,0.1)>,<u_3,(0.9,0.1,0.4)>,\\&
<u_4,(0.4,0.7,0.7)>\}),
(x_2,<u_1,(0.5,0.7,0.7)>,<u_2,(0.5,0.6,0.3)>,\\&<u_3,(0.5,0.6,0.5)>,
<u_4,(0.5,0.8,0.5)>\}),
(x_3,\{<u_1,(0.8,0.6,0.6)>,\\&<u_2,(0.5,0.6,0.7)>,<u_3,(0.7,0.5,0.4)>,
<u_4,(0.3,0.5,0.5)>\})\bigg\} \\
\end{array}
$$
and
$$
\begin{array}{rr}
N_1\tilde{\cap} N_2= &
\bigg\{(x_1,\{<u_1,(0.4,0.6,0.8)>,<u_2,(0.2,0.5,0.8)>,<u_3,(0.3,0.1,0.5)>,\\&
<u_4,(0.4,0.7,0.7)>\}),
(x_2,<u_1,(0.5,0.7,0.8)>,<u_2,(0.3,0.9,0.3)>,\\&<u_3,(0.2,0.6,0.,8)>,
<u_4,(0.4,0.8,0.5)>\}),
(x_3,\{<u_1,(0.7,0.8,0.9)>,\\&<u_2,(0.5,0.9,0.9)>,<u_3,(0.7,0.5,0.8)>,
<u_4,(0.2,0.8,0.6)>\})\bigg\} \\
\end{array}
$$
\end{exmp}
\begin{prop}Let $N_1$, $N_2$  and $N_3$ be any three neutrosophic soft sets. Then,
\begin{enumerate}
\item $ N_1 \widetilde{\cup} N_2= N_2 \widetilde{\cup} N_1$
\item $ N_1 \widetilde{\cap} N_2= N_2 \widetilde{\cap} N_1$
\item $ N_1 \widetilde{\cup}(N_2 \widetilde{\cup} N_3)= (N_1 \widetilde{\cup} N_2) \widetilde{\cup}
N_3$
\item $ N_1 \widetilde{\cap}(N_2 \widetilde{\cap} N_3)= (N_1 \widetilde{\cap} N_2) \widetilde{\cap}
N_3$
\end{enumerate}
\end{prop}
\begin{pf}The proofs can be easily obtained since the t-norm function and s-norm functions are commutative
and associative.
\end{pf}

\subsection{\textbf{Comparision of the Definitions}}

 In this subsection, we compared our definitions of neutrosophic soft with the definitions given
 Maji\cite{maj-13} by inspiring from \cite{cag-14}.

 Let us compare our definitions of neutrosophic soft with the definitions given
 Maji\cite{maj-13} in Table 1.
 $$
\begin{tabular}{|ll|}
  \hline
  \(\textrm{In this paper our approach} \)&
\(\textrm{in  Maji}\)\\
    \hline&  \\
\(N= \{(x, f_N(x)): x\in E\}\)& \(N= \{(x, f_N(x)): x\in A\} \) \\
\(where\)& \(N= \{(x, f_N(x)): x\in A\} \) \\
\(\textrm{E parameter set and}\)& \(A\subseteq E \) \\
\(f_N: E\to N(U)\)& \(f_N: A\to N(U)\} \) \\
     \hline
\end{tabular}
$$

\begin{center}
\footnotesize{\emph{\emph{Table 1}}}
\end{center}
%

Let us compare our complement definitions of neutrosophic soft with
the definitions given
 Maji\cite{maj-13} in Table 2.
$$
\begin{tabular}{|ll|}
  \hline
  \(\textrm{In this paper our approach} \)&
\(\textrm{in  Maji}\)\\
    \hline&  \\
\(N_1^c\)& \(N_1^{\circ}\) \\
\(f_N^c: E\to N(U)\)& \(f_{N_1}^{\circ}: \neg E\to N(U) \) \\
\(T_{f_{N_1^c(x)}}(u)=F_{f_{N_1(x)}}(u)\)& \(T_{f_{N_1^{\circ}(x)}}(u)=F_{f_{N_1(x)}}(u)\) \\
\(I_{f_{N_1^c(x)}}(u)=1-I_{f_{N_1(x)}}(u)\)& \(I_{f_{N_1^{\circ}(x)}}(u)=I_{f_{N_1(x)}}(u)\) \\
\(F_{f_{N_1^c(x)}}(u)=T_{f_{N_1(x)}}(u)\)& \(F_{f_{N_1^{\circ}(x)}}(u)=T_{f_{N_1(x)}}(u)\) \\
     \hline
\end{tabular}
$$
\begin{center}
\footnotesize{\emph{\emph{Table 2}}}
\end{center}

 Let us compare our union definitions of neutrosophic soft with the definitions given
 Maji\cite{maj-13} in Table 2.
$$ \small \begin{tabular}{|ll|}
  \hline
  \(\textrm{\small In this paper our approach} \)&
\(\textrm{in  Maji}\)\\
    \hline&  \\
\(N_3=N_1\tilde{\cup} N_2\)& \(N_3=N_1\acute{{\cup} }N_2 \) \\
\(f_{N_3}: E\to N(U)\)& \(f_{N_3(x)}: A\to N(U) \) \\
\(where\)& \( \) \\
\(T_{f_{N_3(x)}}(u)=s(T_{f_{N_1(x)}}(u),T_{f_{N_2(x)}}(u)) \)&
\(T_{f_{N_3(x)}}(u)=\left \{\begin{array}{ll}
  \displaystyle T_{f_{N_1(x)}}(u) , & x\in A-B
  \\
  \displaystyle T_{f_{N_2(x)}}(u) , &  x\in B- A \\
  \displaystyle max\{T_{f_{N_1(x)}}(u),T_{f_{N_2(x)}}(u)\} , &  x\in A\cap B
\end{array}\right. \) \\
\(I_{f_{N_3(x)}}(u)=t(I_{f_{N_1(x)}}(u),I_{f_{N_2(x)}}(u))\)&
\(I_{f_{N_3(x)}}(u)=\left \{\begin{array}{ll}
  \displaystyle I_{f_{N_1(x)}} (u), & x\in A-B
  \\
  \displaystyle I_{f_{N_2(x)}}(u) , &  x\in B- A \\
  \displaystyle \frac{(I_{f_{N_1(x)}}(u),I_{f_{N_2(x)}}(u))}{2} ,\,\,\,\,\,\,\,\,\,\,\,\,\, &  x\in A\cap B
\end{array}\right.\) \\
\(F_{f_{N_3(x)}}(u)=t(F_{f_{N_1(x)}}(u),F_{f_{N_2(x)}}(u))\)&
\(F_{f_{N_3(x)}}(u)=\left \{\begin{array}{ll}
  \displaystyle F_{f_{N_1(x)}}(u) , & x\in A-B
  \\
  \displaystyle F_{f_{N_2(x)}}(u) , &  x\in B- A \\
  \displaystyle min\{I_{f_{N_1(x)}}(u),I_{f_{N_2(x)}}(u)\} , &  x\in A\cap B
\end{array}\right. \) \\
     \hline
\end{tabular}
$$
\begin{center}
\footnotesize{\emph{\emph{Table 2}}}
\end{center}

 Let us compare our intersection definitions of neutrosophic soft with the definitions given
 Maji\cite{maj-13} in Table 2.
$$
\small  \begin{tabular}{|ll|}
  \hline
  \(\textrm{\small In this paper our approach} \)&
\(\textrm{in  Maji}\)\\
    \hline&  \\
\(N_3=N_1\tilde{\cap} N_2\)& \(N_3=N_1\acute{{\cap} }N_2 \) \\
\(f_{N_3}: E\to N(U)\)& \(f_{N_3(x)}: A\to N(U) \) \\
\(where\)& \( \) \\
\(T_{f_{N_3(x)}}(u)=t(T_{f_{N_1(x)}}(u),T_{f_{N_2(x)}}(u)) \)&
\( T_{f_{N_3(x)}}(u)=min\{T_{f_{N_1(x)}}(u),T_{f_{N_2(x)}}(u)\} \) \\
\(I_{f_{N_3(x)}}(u)=s(I_{f_{N_1(x)}}(u),I_{f_{N_2(x)}}(u))\)&
\(I_{f_{N_3(x)}}(u)=\frac{(I_{f_{N_1(x)}}(u),I_{f_{N_2(x)}}(u))}{2}\) \\
\(F_{f_{N_3(x)}}(u)=s(F_{f_{N_1(x)}}(u),F_{f_{N_2(x)}}(u))\)&
\(F_{f_{N_3(x)}}(u)=max\{F_{f_{N_1(x)}}(u),F_{f_{N_2(x)}}(u)\} \) \\
     \hline
\end{tabular}
$$
\begin{center}

\footnotesize{\emph{\emph{Table 3}}}
\end{center}
\section{Neutrosophic Soft Matrices}
In this section, we presented neutrosophic soft matrices which are
representative of the neutrosophic soft sets. The matrix is useful
for storing an neutrosophic soft set in computer memory which are
very useful and applicable. Some of it is quoted from
\cite{cag-09sm, cag-12, bas-13}.

This section are an attempt to extend the concept of soft
matrices\cite{cag-09sm}, fuzzy soft matrices\cite{cag-12},
intuitionistic fuzzy soft matrices\cite{bas-13}.
\begin{defn}
Let $N$ be an neutrosophic soft set over $N(U)$. Then a subset of
$N(U)\times E$ is uniquely defined by

$R_N=\{({f_N(x)},x):x \in E, f_N(x)\in N(U)\}$ which is called a
relation form of $(N,E)$. The characteristic function of $R_N$ is
written by
$$
\Theta_{R_N}:N(U)\times E\to [0,1]\times[0,1]\times[0,1],
\Theta_{R_N}(u,x)=(T_{f_N(x)}(u),I_{f_N(x)}(u),F_{f_N(x)}(u))
$$
where $T_{f_N(x)}(u)$, $I_{f_N(x)}(u)$ and $F_{f_N(x)}(u)$ is the
truth-membership, indeterminacy-membership and falsity-membership of
$u\in U$, respectively.
\end{defn}

\begin{defn}
Let $U=\{u_1,u_2,\ldots,u_m\}$, $E=\{x_1,x_2,\ldots,x_n\}$ and $N$
be an neutrosophic soft set over $N(U)$. Then
$$
\begin{tabular}{c|cccc}
$R_N$ & $ f_N(x_1)$         & $ f_N(x_2)$     & $\cdots$ & $ f_N(x_n)$\\
\hline
$u_1$     & $\Theta_{R_N}(u_1,x_1)$ & $\Theta_{R_N}(u_1,x_2)$ & $\cdots$ & $\Theta_{R_N}(u_1,x_n)$\\
$u_2$     & $\Theta_{R_N}(u_2,x_1)$ & $\Theta_{R_N}(u_2,x_2)$ & $\cdots$ & $\Theta_{R_N}(u_2,x_n)$\\
$\vdots$  & $\vdots$              & $\vdots$              & $\ddots$ & $\vdots$ \\
$u_m$     & $\Theta_{R_N}(u_m,x_1)$ & $\Theta_{R_N}(u_m,x_2)$ &
$\cdots$ & $\Theta_{R_N}(u_m,x_n)$
\end{tabular}
$$
\end{defn}
If $a_{ij}=\Theta_{R_N}(u_i,x_j)$, we can define a matrix
$$
[a_{ij}]=\begin{bmatrix}
a_{11} & a_{12} & \cdots & a_{1n} \\
a_{21} & a_{22} & \cdots & a_{2n} \\
\vdots & \vdots & \ddots & \vdots \\
a_{m1} & a_{m2} & \cdots & a_{mn}
\end{bmatrix}
$$
such that
$a_{ij}=(T_{f_N(x_j)}(u_i),I_{f_N(x_j)}(u_i),F_{f_N(x_j)}(u_i))=(T_{ij}^a,I_{ij}^a,F_{ij}^a)$,
which is called an $m\times n$ neutrosophic soft matrix (or namely
NS-matrix) of the neutrosophic soft set $N$ over $N(U)$.

According to this definition, an a neutrosophic soft set $N$ is
uniquely characterized by matrix $[a_{ij}]_{m\times n}$. Therefore,
we shall identify any neutrosophic soft set with its soft NS-matrix
and use these two concepts as interchangeable. The set of all
$m\times n$ NS-matrix over $N(U)$ will be denoted by
$\widetilde{N}_{m\times n}$. From now on we shall delete th
subscripts $m\times n$ of $[a_{ij}]_{m\times n}$, we use $[a_{ij}]$
instead of $[a_{ij}]_{m\times n}$, since
$[a_{ij}]\in\widetilde{N}_{m\times n}$ means that $[a_{ij}]$ is an
$m\times n$ NS-matrix for $i=1,2,\ldots,m$ and $j=1,2,\ldots,n$.

\begin{exmp}\label{x1}
Let $U=\{u_1, u_2, u_3\}$, $E=\{x_1, x_2, x_3\}$. $N_1$ be a
neutrosophic soft sets over neutrosophic as
$$
\begin{array}{rr}
N&=
\bigg\{(x_1,\{<u_1,(0.7,0.6,0.7)>,<u_2,(0.4,0.2,0.8)>,<u_3,(0.9,0.1,0.5)>
\}),\\&
(x_2,<u_1,(0.5,0.7,0.8)>,<u_2,(0.5,0.9,0.3)>,<u_3,(0.5,0.6,0.,8)>
\}),\\&
(x_3,\{<u_1,(0.8,0.6,0.9)>,<u_2,(0.5,0.9,0.9)>,<u_3,(0.7,0.5,0.4)>
\})\bigg\} \\
\end{array}
$$
Then, the NS-matrix  $[a_{ij}]$ is written by
$$
[a_{ij}]=\begin{bmatrix}
(0.7,0.6,0.7) & (0.5,0.7,0.8) & (0.8,0.6,0.9) \\
(0.4,0.2,0.8) & (0.5,0.9,0.3) & (0.5,0.9,0.9) \\
(0.9,0.1,0.5) & (0.5,0.6,0.8) & (0.7,0.5,0.4) \\
\end{bmatrix}
$$
\end{exmp}

\begin{defn} A neutrosophic soft matrix of order $1\times n$ i.e., with a single
row is called a row-neutrosophic soft matrix. Physically, a
row-neutrosophic soft matrix formally corresponds to an neutrosophic
soft set whose universal set contains only one object.
\end{defn}
\begin{exmp}
Let $U=\{u_1\}$, $E=\{x_1, x_2, x_3\}$. $N_1$ be a neutrosophic soft
sets over neutrosophic as
$$
\begin{array}{rr}
N= & \bigg\{(x_1,\{<u_1,(0.7,0.6,0.7)>\}),
(x_2,<u_1,(0.5,0.7,0.8)>\}),\\&
(x_3,\{<u_1,(0.8,0.6,0.9)>\})\bigg\} \\
\end{array}
$$
Then, the NS-matrix  $[a_{ij}]$ is written by
$$
[a_{ij}]=\begin{bmatrix}
(0.7,0.6,0.7) & (0.5,0.7,0.8) & (0.8,0.6,0.9) \\
\end{bmatrix}.
$$

\end{exmp}

\begin{defn} A neutrosophic soft matrix of order $m\times 1$ i.e., with a
single column is called a column-neutrosophic soft matrix.
Physically, a column-neutrosophic soft matrix formally corresponds
to an neutrosophic soft set whose parameter set contains only one
parameter.
\end{defn}
\begin{exmp}
Let $U=\{u_1, u_2, u_3,u_4\}$, $E=\{x_1\}$. $N_1$ be a neutrosophic
soft sets over neutrosophic as
$$
\begin{array}{rr}
N= &
\bigg\{(x_1,\{<u_1,(0.7,0.6,0.7)>,<u_2,(0.4,0.2,0.8)>,<u_3,(0.9,0.1,0.5)>,\\&
<u_4,(0.4,0.7,0.7)>\})\bigg\} \\
\end{array}
$$
Then, the NS-matrix  $[a_{ij}]$ is written by
$$
[a_{ij}]=\begin{bmatrix}
(0.7,0.6,0.7) \\
(0.4,0.2,0.8)  \\
(0.9,0.1,0.5)  \\
(0.4,0.7,0.7)  \\
\end{bmatrix}.
$$
\end{exmp}

\begin{defn} A neutrosophic soft matrix of order $m \times n$ is said to be a
square neutrosophic soft matrix if $m = n$ i.e., the number of rows
and the number of columns are equal. That means a
square-neutrosophic soft matrix is formally equal to an neutrosophic
soft set having the same number of objects and parameters.
\end{defn}
\begin{exmp}
Consider the Example \ref{x1}. Here since the NS-matrix contains
three rows and three columns, so it is a square-neutrosophic soft
matrix.
\end{exmp}

\begin{defn} A square neutrosophic soft matrix of order $m\times n$ is said to
be a diagonal-neutrosophic soft matrix if all of its non-diagonal
elements are $ (0,0,1).$
\end{defn}
\begin{exmp}
Let $U=\{u_1, u_2, u_3,u_4\}$, $E=\{x_1, x_2, x_3\}$. $N_1$ be a
neutrosophic soft sets over neutrosophic as
$$
\begin{array}{rr}
N= &
\bigg\{(x_1,\{<u_1,(0.7,0.6,0.7)>,<u_2,(0.0,1.0,1.0)>,<u_3,(0.0,1.0,1.0)>\}),\\&
(x_2,<u_1,(0.0,1.0,1.0)>,<u_2,(0.0,1.0,1.0)
>,<u_3,(0.0,1.0,1.0)>\}),\\&
(x_3,\{<u_1,(0.0,1.0,1.0)>,<u_2,(0.0,1.0,1.0)>,<u_3,(0.7,0.5,0.4)>\})\bigg\} \\
\end{array}
$$
Then, the NS-matrix  $[a_{ij}]$ is written by
$$
[a_{ij}]=\begin{bmatrix}
(0.7,0.6,0.7) & (0.0,1.0,1.0) & (0.0,1.0,1.0) \\
(0.0,1.0,1.0) & (0.0,1.0,1.0) & (0.0,1.0,1.0) \\
(0.0,1.0,1.0) & (0.0,1.0,1.0) & (0.7,0.5,0.4) \\
\end{bmatrix}.
$$
\end{exmp}
\begin{defn} The transpose of a square neutrosophic soft matrix $[a_{ij}]$ of
order $m \times n$ is another square neutrosophic soft matrix of
order $n \times m$ obtained from $[a_{ij}]$ by interchanging its
rows and columns. It is denoted by $[a_{ij}^T]$. Therefore the
neutrosophic soft set associated with $[a_{ij}^T]$ becomes a new
neutrosophic soft set over the same universe and over the same set
of parameters.
\end{defn}
\begin{exmp}
Consider the Example \ref{x1}. If the NS-matrix  $[a_{ij}]$ is
written by
$$
[a_{ij}]=\begin{bmatrix}
(0.7,0.6,0.7) & (0.5,0.7,0.8) & (0.8,0.6,0.9) \\
(0.4,0.2,0.8) & (0.5,0.9,0.3) & (0.5,0.9,0.9) \\
(0.9,0.1,0.5) & (0.5,0.6,0.8) & (0.7,0.5,0.4) \\
\end{bmatrix}.
$$
then, its transpose neutrosophic soft matrix as;
$$
[a_{ij}]=\begin{bmatrix}
(0.7,0.6,0.7) & (0.4,0.2,0.8) &  (0.9,0.1,0.5) \\
(0.5,0.7,0.8) & (0.5,0.9,0.3) &  (0.5,0.6,0.8) \\
(0.8,0.6,0.9) & (0.5,0.9,0.9) &  (0.7,0.5,0.4) \\
\end{bmatrix}.
$$
\end{exmp}

\begin{defn} A square neutrosophic soft matrix $[a_{ij}]$ of order $n\times n$ is said
to be a symmetric neutrosophic soft matrix, if its transpose be
equal to it, i.e., if $[a_{ij}^T]=[a_{ij}]$. Hence the neutrosophic
soft matrix $[a_{ij}]$) is symmetric, if $[a_{ij}]$= $[a_{ji}]$
$\forall i, j.$
\end{defn}
\begin{exmp}
Let $U=\{u_1, u_2, u_3\}$, $E=\{x_1, x_2, x_3\}$. $N_1$ be a
neutrosophic soft sets as
$$
\begin{array}{rr}
N= &
\bigg\{(x_1,\{<u_1,(0.7,0.6,0.7)>,<u_2,(0.4,0.2,0.8)>,<u_3,(0.9,0.1,0.5)>\}),\\&
(x_2,<u_1,(0.4,0.2,0.8)>,<u_2,(0.5,0.9,0.3)>,<u_3,(0.5,0.9,0.9)>\}),\\&
(x_3,\{<u_1,(0.9,0.1,0.5)>,<u_2,(0.5,0.9,0.9)>,<u_3,(0.7,0.5,0.4)>\})\bigg\} \\
\end{array}
$$
Then, the symmetric neutrosophic matrix  $[a_{ij}]$ is written by
$$
[a_{ij}]=\begin{bmatrix}
(0.7,0.6,0.7) & (0.4,0.2,0.8) & (0.9,0.1,0.5) \\
(0.4,0.2,0.8) & (0.5,0.9,0.3) & (0.5,0.9,0.9)\\
(0.9,0.1,0.5) & (0.5,0.9,0.9) & (0.7,0.5,0.4) \\
\end{bmatrix}.
$$
\end{exmp}
\begin{defn}
Let $[a_{ij}]\in \widetilde{N}_{m\times n}$. Then $[a_{ij}]$ is
called
\begin{enumerate}[i.]
\item A zero NS-matrix, denoted by $[\tilde 0]$, if $a_{ij}=(0,1,1)$ for all $i$ and $j$.
\item A universal NS-matrix, denoted by $[\tilde 1]$, if $a_{ij}=(1,0,0)$
for all $i$ and $j$.
\end{enumerate}
\end{defn}

\begin{exmp}
Let $U=\{u_1, u_2, u_3\}$, $E=\{x_1, x_2, x_3\}$. Then, a zero
NS-matrix $[a_{ij}]$ is written by
$$
[a_{ij}]=\begin{bmatrix}
(0,1,1) & (0,1,1) & (0,1,1) \\
(0,1,1) & (0,1,1) & (0,1,1) \\
(0,1,1) & (0,1,1) & (0,1,1) \\
\end{bmatrix}.
$$
and a universal NS-matrix $[a_{ij}]$ is written by
$$
[a_{ij}]=\begin{bmatrix}
(1,0,0) & (1,0,0) & (1,0,0) \\
(1,0,0) & (1,0,0) & (1,0,0) \\
(1,0,0) & (1,0,0) & (1,0,0) \\
\end{bmatrix}.
$$
\end{exmp}
\begin{defn}
Let $[a_{ij}],[b_{ij}]\in \widetilde{N}_{m\times n}$. Then
\begin{enumerate}[i.]
\item $[a_{ij}]$ is an NS-submatrix of $[b_{ij}]$, denoted,
$[a_{ij}]\tilde \subseteq[b_{ij}]$, if $T_{ij}^b\geq T_{ij}^a$,
$I_{ij}^a\geq I_{ij}^b$ and $F_{ij}^a\geq F_{ij}^b$, for all $i$ and
$j$.
\item $[a_{ij}]$ is a proper NS-submatrix of $[b_{ij}]$, denoted,
$[a_{ij}]\tilde \subset [b_{ij}]$, if $T_{ij}^a\geq T_{ij}^b$,
$I_{ij}^a\leq I_{ij}^b$ and $F_{ij}^a\leq F_{ij}^b$ for at least
$T_{ij}^a>T_{ij}^b$ and $I_{ij}^a< I_{ij}^b$ and $F_{ij}^a<
F_{ij}^b$ for all $i$ and $j$.
\item $[a_{ij}]$ and $[b_{ij}]$ are IFS equal matrices, denoted by
$[a_{ij}]=[b_{ij}]$, if $a_{ij}=b_{ij}$ for all $i$ and $j$.
\end{enumerate}
\end{defn}
\begin{defn}
Let $[a_{ij}],[b_{ij}]\in \widetilde{N}_{m\times n}$. Then
\begin{enumerate}[i.]
\item Union of $[a_{ij}]$ and  $[b_{ij}]$, denoted, $[a_{ij}]\tilde\cup[b_{ij}]$,
if $c_{ij}=(T^c_{ij},I^c_{ij},F^c_{ij})$, where
$T_{ij}^c=\max\{T_{ij}^a,T_{ij}^b\}$,
$I_{ij}^c=\min\{I_{ij}^a,I_{ij}^b\}$ and
$F_{ij}^c=\min\{F_{ij}^a,F_{ij}^b\}$ for all $i$ and $j$.
\item Intersection of $[a_{ij}]$ and  $[b_{ij}]$, denoted, $[a_{ij}]\tilde\cap[b_{ij}]$,
if $c_{ij}=(T^c_{ij},I^c_{ij},F^c_{ij})$, where
$T_{ij}^c=\min\{T_{ij}^a,T_{ij}^b\}$,
$I_{ij}^c=\max\{I_{ij}^a,I_{ij}^b\}$ and
$F_{ij}^c=\max\{F_{ij}^a,F_{ij}^b\}$ for all $i$ and $j$.
\item Complement of $[a_{ij}]$, denoted by $[a_{ij}]^c$, if $c_{ij}=(F_{ij}^a,1-I^a_{ij},T^a_{ij})$
for all $i$ and $j$.
\end{enumerate}
\end{defn}

\begin{exmp}Consider the Example \ref{1}. Then,
$$
[a_{ij}]\tilde\cup[b_{ij}]=\begin{bmatrix}
(0.7,0.5,0.7) & (0.5,0.7,0.7) & (0.8,0.6,0.6) \\
 (0.4,0.2,0.1) & (0.5,0.6,0.3)  & (0.5,0.6,0.7) \\
 (0.9,0.1,0.4) & (0.5,0.6,0.5) & (0.7,0.5,0.4) \\
  (0.4,0.7,0.7) & (0.5,0.8,0.5) &  (0.3,0.5,0.5) \\
\end{bmatrix},
$$

$$
[a_{ij}]\tilde\cap[b_{ij}]=\begin{bmatrix}
(0.4,0.6,0.8) & (0.5,0.7,0.8) & (0.7,0.8,0.9) \\
 (0.2,0.5,0.8) & (0.3,0.9,0.3)  & (0.5,0.9,0.9) \\
(0.3,0.1,0.5) & (0.2,0.6,0.8) & (0.7,0.5,0.8) \\
  (0.4,0.7,0.7) & (0.4,0.8,0.5) &  (0.2,0.8,0.6) \\
\end{bmatrix}.
$$
and

$$
[a_{ij}]^c=\begin{bmatrix}
(0.8,0.5,0.4) & (0.7,0.3,0.5) & (0.6,0.2,0.7) \\
(0.1,0.5,0.2) & (0.3,0.4,0.3) & (0.7,0.4,0.5) \\
(0.4,0.9,0.3) & (0.5,0.4,0.2) & (0.8,0.5,0.7) \\
(0.7,0.3,0.4) & (0.5,0.5,0.4) & (0.5,0.2,0.2) \\
\end{bmatrix}.
$$

\end{exmp}
\begin{defn}
Let $[a_{ij}],[b_{ij}]\in \widetilde{N}_{m\times n}$. Then
$[a_{ij}]$ and $[b_{ij}]$ are disjoint, if
$[a_{ij}]\tilde\cap[b_{ij}]=[\tilde 0]$ for all $i$ and $j$.
\end{defn}
\begin{prop}
Let $[a_{ij}]\in \widetilde{N}_{m\times n}$. Then
\begin{enumerate}[i.]
\item $\big([a_{ij}]^c\big)^c=[a_{ij}]$
\item $[\tilde 0]^c=[\tilde 1]$.
\end{enumerate}
\end{prop}
\begin{prop}
Let $[a_{ij}],[b_{ij}]\in \widetilde{N}_{m\times n}$. Then
\begin{enumerate}[i.]
\item $[a_{ij}]\subseteq[\tilde 1]$
\item $[\tilde 0]\tilde\subseteq[a_{ij}]$
\item $[a_{ij}]\tilde\subseteq[a_{ij}]$
\item $[a_{ij}]\tilde\subseteq[b_{ij}]$ and $[b_{ij}]\tilde\subseteq[c_{ij}]$
$\Rightarrow$ $[a_{ij}]\tilde\subseteq[c_{ij}]$
\end{enumerate}
\end{prop}
\begin{prop}
Let $[a_{ij}],[b_{ij}],[c_{ij}]\in \widetilde{N}_{m\times n}$. Then
\begin{enumerate}[i.]
\item $[a_{ij}]=[b_{ij}]$ and $[b_{ij}]=[c_{ij}]$
$\Leftrightarrow$ $[a_{ij}]=[c_{ij}]$
\item $[a_{ij}]\tilde\subseteq[b_{ij}]$ and $[b_{ij}]\tilde\subseteq[a_{ij}]$
$\Leftrightarrow$ $[a_{ij}]=[b_{ij}]$
\end{enumerate}
\end{prop}
\begin{prop}
Let $[a_{ij}],[b_{ij}],[c_{ij}]\in \widetilde{N}_{m\times n}$. Then
\begin{enumerate}[i.]
\item $[a_{ij}]\tilde\cup[a_{ij}]=[a_{ij}]$
\item $[a_{ij}]\tilde\cup[\tilde 0]=[a_{ij}]$
\item $[a_{ij}]\tilde\cup[\tilde 1]=[\tilde 1]$
\item $[a_{ij}]\tilde\cup[b_{ij}]=[b_{ij}]\tilde\cup[a_{ij}]$
\item $([a_{ij}]\tilde\cup[b_{ij}])\tilde\cup[c_{ij}]=
[a_{ij}]\tilde\cup([b_{ij}]\tilde\cup[c_{ij}])$
\end{enumerate}
\end{prop}
\begin{prop}
Let $[a_{ij}],[b_{ij}],[c_{ij}]\in\widetilde{N}_{m\times n}$. Then
\begin{enumerate}[i.]
\item $[a_{ij}]\tilde\cap[a_{ij}]=[a_{ij}]$
\item $[a_{ij}]\tilde\cap[\tilde 0]=[\tilde 0]$
\item $[a_{ij}]\tilde\cap[\tilde 1]=[a_{ij}]$
\item $[a_{ij}]\tilde\cap[b_{ij}]=[b_{ij}]\tilde\cap[a_{ij}]$
\item $([a_{ij}]\tilde\cap[b_{ij}])\tilde\cap[c_{ij}]=
[a_{ij}]\tilde\cap([b_{ij}]\tilde\cap[c_{ij}])$
\end{enumerate}
\end{prop}
\begin{prop}
Let $[a_{ij}],[b_{ij}]\in \widetilde{N}_{m\times n}$. Then De
Morgan's laws are valid
\begin{enumerate}[i.]
\item $([a_{ij}]\tilde\cup[b_{ij}])^c=[a_{ij}]^c\tilde\cap[b_{ij}]^c$
\item $([a_{ij}]\tilde\cap[b_{ij}])^c=[a_{ij}]^c\tilde\cup[b_{ij}]^c$
\end{enumerate}
\end{prop}
\begin{pf}
\emph{i.}
\begin{eqnarray*}
([a_{ij}]\tilde\cup[b_{ij}])^c & = & ([(T_{ij}^a,I_{ij}^a,F_{ij}^a)]\tilde\cup[(T_{ij}^b,I_{ij}^b,F_{ij}^b)])^c\\
                               & = & [(\max\{T_{ij}^a,T_{ij}^b\},\min\{I_{ij}^a,I_{ij}^b\},\min\{F_{ij}^a,F_{ij}^b\})]^c\\
                               & = & [(\min\{F_{ij}^a,F_{ij}^b\},\max\{1-I_{ij}^a,1-I_{ij}^b\},\max\{T_{ij}^a,T_{ij}^b\}))] \\
                               & = & [(F_{ij}^a,I_{ij}^a,T_{ij}^a)]\tilde\cap[(F_{ij}^b,I_{ij}^b,T_{ij}^b)]\\
                               & = & [a_{ij}]^c\tilde\cap[b_{ij}]^c
\end{eqnarray*}
\emph{i.}
\begin{eqnarray*}
([a_{ij}]\tilde\cap[b_{ij}])^c & = & ([(T_{ij}^a,I_{ij}^a,F_{ij}^a)]\tilde\cap[(T_{ij}^b,I_{ij}^b,F_{ij}^b)])^c\\
                               & = & [(\min\{T_{ij}^a,T_{ij}^b\},\max\{I_{ij}^a,I_{ij}^b\},\max\{F_{ij}^a,F_{ij}^b\})]^c\\
                               & = & [(\max\{F_{ij}^a,F_{ij}^b\},\min\{1-I_{ij}^a,1-I_{ij}^b\},\min\{T_{ij}^a,T_{ij}^b\}))] \\
                               & = & [(F_{ij}^a,I_{ij}^a,T_{ij}^a)]\tilde\cup[(F_{ij}^b,I_{ij}^b,T_{ij}^b)]\\
                               & = & [a_{ij}]^c\tilde\cup[b_{ij}]^c
\end{eqnarray*}
\end{pf}
\begin{prop}
Let $[a_{ij}],[b_{ij}],[c_{ij}]\in\widetilde{N}_{m\times n}$. Then
\begin{enumerate}[i.]
\item $[a_{ij}]\tilde\cap([b_{ij}]\tilde\cup[c_{ij}])=
([a_{ij}]\tilde\cap([b_{ij}])\tilde\cup([a_{ij}]\tilde\cap[c_{ij}])$
\item $[a_{ij}]\tilde\cup([b_{ij}]\tilde\cap[c_{ij}])=
([a_{ij}]\tilde\cup([b_{ij}])\tilde\cap([a_{ij}]\tilde\cup[c_{ij}])$
\end{enumerate}
\end{prop}

\section{Products of NS-Matrices}
In this section, we define two special products of NS-matrices to
construct soft decision making methods.

\begin{defn}
Let $[a_{ij}],[b_{ik}]\in\widetilde{N}_{m\times n}$. Then
\emph{And}-product of $[a_{ij}]$ and $[b_{ij}]$  is defined by
$$
\wedge:\widetilde{N}_{m\times n}\times\widetilde{N}_{m\times n}
\to\widetilde{N}_{m\times n^2}\qquad
[a_{ij}]\wedge[b_{ik}]=[c_{ip}]=(T_{ip}^c,I_{ip}^c,F_{ip}^c)
$$
where

$T_{ip}^c=t(T_{ij}^a,T_{jk}^b) $, $ I_{ip}^c=s(I_{ij}^a,I_{jk}^b) $
and $F_{ip}^c=s(F_{ij}^a,F_{jk}^b) $ such that $p=n(j-1)+k$
\end{defn}
\begin{defn}
Let $[a_{ij}],[b_{ik}]\in\widetilde{N}_{m\times n}$. Then
\emph{And}-product of $[a_{ij}]$ and $[b_{ij}]$  is defined by
$$
\vee:\widetilde{N}_{m\times n}\times\widetilde{N}_{m\times n}
\to\widetilde{N}_{m\times n^2}\qquad
[a_{ij}]\vee[b_{ik}]=[c_{ip}]=(T_{ip}^c,I_{ip}^c,F_{ip}^c)
$$
where

$T_{ip}^c=s(T_{ij}^a,T_{jk}^b) $, $ I_{ip}^c=t(I_{ij}^a,I_{jk}^b) $
and $F_{ip}^c=t(F_{ij}^a,F_{jk}^b) $ such that $p=n(j-1)+k$
\end{defn}
\begin{exmp}
Assume that $[a_{ij}],[b_{ik}]\in\widetilde{N}_{3\times 2}$ are
given as follows
$$
[a_{ij}]=\begin{bmatrix}
(1.0,0.1,0.1) & (1.0,0.4,0.1) &  \\
(1.0,0.2,0.1) & (1.0,0.1,0.1) & \\
(1.0,0.8,0.1) & (1.0,0.7,0.1) &
\end{bmatrix}
$$
$$
[b_{ij}]=\begin{bmatrix}
(1.0,0.7,0.1) & (1.0,0.1,0.1)  \\
(1.0,0.5,0.1) & (1.0,0.2,0.1)  \\
(1.0,0.5,0.1) & (1.0,0.5,0.1)
\end{bmatrix}
$$
$$
[a_{ij}]\wedge[b_{ij}]=\begin{bmatrix}
(1.0,0.7,0.1) & (1.0,0.1,0.1) & (1.0,0.7,0.1)&(1.0,0.4,0.1)  \\
(1.0,0.5,0.1) & (1.0,0.2,0.1) & (1.0,0.5,0.1)&(1.0,0.2,0.1) \\
(1.0,0.8,0.1) & (1.0,0.8,0.1) & (1.0,0.7,0.1)&(1.0,0.7,0.1)
\end{bmatrix}
$$
$$
[a_{ij}]\vee[b_{ij}]=\begin{bmatrix}
(1.0,0.1,0.1) & (1.0,0.1,0.1) & (1.0,0.4,0.1)&(1.0,0.1,0.1)  \\
(1.0,0.2,0.1) & (1.0,0.2,0.1) & (1.0,0.1,0.1)&(1.0,0.1,0.1) \\
(1.0,0.8,0.1) & (1.0,0.8,0.1) & (1.0,0.5,0.1)&(1.0,0.5,0.1)
\end{bmatrix}
$$
\end{exmp}
\begin{prop}
Let $[a_{ij}],[b_{ij}],[c_{ij}]\in\widetilde{N}_{m\times n}$. Then
the De morgan's types of results are true.
\begin{enumerate}[i.]
\item $([a_{ij}]\vee[b_{ij}])^c=[a_{ij}]^c\wedge[b_{ij}]^c$
\item $([a_{ij}]\wedge[b_{ij}])^c=[a_{ij}]^c\vee[b_{ij}]^c$
\end{enumerate}
\end{prop}
\section{Decision making problem using and-product of neutrosophic soft
matrices}

\begin{defn}\label{SMm}\cite{cag-09sm}
Let $[(\mu_{ip},\nu_{ip},w_{ip})]\in NSM_{m\times n^2}$, $I_k=\{p:
(\mu_{ip},\nu_{ip},w_{ip})\neq 0,$ for some $1\leq i\leq m,$
$(k-1)n<p\leqslant kn$ for all $k\in I=\{1,2,...,n\}$. Then
$NS$-max-min decision function, denoted $D_{mMM}$, is defined as
follows
$$
D_{mMM}:NSM_{m\times n^2}\to NSM_ {m\times 1},$$
$$
D_{mMM}=[(\mu_{ip},\nu_{ip},w_{ip})]=[d_{i1}]=[(\max_{k}\{\acute{\mu_{ipk}}\},\{\acute{\nu_{ipk}}\},\min_{k}\{\acute{w_{ipk}}\})]
$$
where
$$
\acute{\mu_{ipk}}= \left\{
  \begin{array}{ll}
    \max_{p\in  I_k}\{\mu_{ipk}\},&\hbox{if }  I_k\neq \emptyset,\\
    0,&\hbox{if }  I_k=\emptyset.
  \end{array}
\right.
$$

$$
\acute{\nu_{ipk}}= \left\{
  \begin{array}{ll}
    \min_{p\in  I_k}\{\nu_{ipk}\},&\hbox{if }  I_k\neq \emptyset,\\
    0,&\hbox{if }  I_k=\emptyset.
  \end{array}
\right.
$$

$$
\acute{w_{ipk}}= \left\{
  \begin{array}{ll}
    \min_{p\in  I_k}\{w_{ipk}\},&\hbox{if }  I_k\neq \emptyset,\\
    0,&\hbox{if }  I_k=\emptyset.
  \end{array}
\right.
$$
The one column $fs$-matrix $Mm[c_{ip}]$ is called max-min decision
$fs$-matrix.
\end{defn}

\begin{defn} Let $U = \{u_1,u_2,u_3,u_m\}$ be the universe
and $D_{mMM}(\mu_{ip},\nu_{ip},w_{ip})= [d_{i1}]$. Then the set
defined by

$$opt_{[d_{i1}]}^m (U)=\{u_i/d_i : u_i\in U, d_i= max \{s_i\}\},$$
where $s_i=\mu_{ip}-\nu_{ip},w_{ip}, d_{i1}\neq 0$ which is called
an optimum fuzzy set on U.
\end{defn}
\textbf{Algorithm }

The algorithm for the solution is given below

\textbf{Step 1}: Choose feasible subset of the set of parameters.

\textbf{Step 2}: Construct the neutrosophic matrices for each
parameter.

\textbf{Step 3:} Choose a product of the neutrosophic matrices ,

\textbf{Step 4:} Find the method min-max-max decision N-matrices.

\textbf{Step 5:} Find an optimum fuzzy set on U.

\begin{rem} We can also define NS-matrices max-min-min decision making
methods. One of them may be more useful than the others according to
the type of problem.
\end{rem}

\textbf{Case study:} Assume that , a car dealer stores three
different types of cars $U=\{u_1,u_2, u_3\}$ which may be
characterize by the set of parameters $E=\{e_1 ,  e_2\}$ where $e_1$
stands for costly , $e_2$ stands for fuel effciency. Then we
consider the following example. Suppose a couple Mr. X and  Mrs. X
come to the dealer to buy a car before Dugra Puja. If each partner
has to consider his/her own set of parameters, then we select the
car on the basis of partner's parameters by using NS-matrices
min-max-max decision making as follow.

\textbf{Step 1}: First Mr.X and Mrs.X have to chose the sets of
their parameters $A= \{e_1 ,  e_2\}$ and $B =\{e_1 ,  e_2\}$,
respectively.

\textbf{Step 2}:Then we construct the NS-matrices $[a_{ij}]$ and
$[b_{ij}]$ according to their set of parameters A and B,
respectively, as follow:
 $$
[a_{ij}]=\begin{bmatrix}
(1.0,0.1,0.1) & (1.0,0.4,0.1) &  \\
(1.0,0.2,0.1) & (1.0,0.1,0.1) & \\
(1.0,0.8,0.1) & (1.0,0.7,0.1) &
\end{bmatrix}
$$
and
$$
[b_{ij}]=\begin{bmatrix}
(1.0,0.7,0.1) & (1.0,0.1,0.1)  \\
(1.0,0.5,0.1) & (1.0,0.2,0.1)  \\
(1.0,0.5,0.1) & (1.0,0.5,0.1)
\end{bmatrix}
$$

\textbf{Step 3:}Now ,we can find the And-product of the NS-matrices
 $[a_{ij}]$ and
$[b_{ij}]$ as follow:
$$
[a_{ij}]\wedge[b_{ij}]=\begin{bmatrix}
(1.0,0.7,0.1) & (1.0,0.1,0.1) & (1.0,0.7,0.1)&(1.0,0.4,0.1)  \\
(1.0,0.5,0.1) & (1.0,0.2,0.1) & (1.0,0.5,0.1)&(1.0,0.2,0.1) \\
(1.0,0.8,0.1) & (1.0,0.8,0.1) & (1.0,0.7,0.1)&(1.0,0.7,0.1)
\end{bmatrix}
$$

\textbf{Step 4: } Now,we calculate; for $i=\{1,2,3\}$
 $$
[d_{i1}]=\begin{bmatrix}
(\mu_{11},\nu_{11},w_{11})&   \\
(\mu_{21},\nu_{21},w_{21})& \\
(\mu_{31},\nu_{31},w_{31}) &
\end{bmatrix}
$$
To demonstrate, let us  find  $d_{21}$ for $i=2$. Since $i=2$ and $k
\in\{1,2\}$ so $d_{21}=(\mu_{21},\nu_{21},w_{21})$.

Let $t_{2k} =\{ t_{21} , t_{22}   \}$, where
$t_{2k}=(\mu_{2p},\nu_{2p},w_{2p})$ then,

we have to find $t_{2k}$ for all $k \in\{1,2\}$. First to find
$t_{21}$, $I_1=\{ p: 0< p \leq 2\}$ for $k=1$ and $n=2$.

We have $t_{21} =( min\{ \mu_{2p}\}, max\{\nu_{2p} \}, max\{w_{2p}
\})$,

here $p =1,2$ $(min\{ \mu_{21},\mu_{22}\}, max\{\nu_{21},\nu_{22}\},
max\{w_{21},w_{22} \})\\ =(min\{ 1,1\}, max\{0.5,0.2\},
max\{0.1,0.1\})=(1 ,0.5 , 0.1)$
 and

 $t_{22} =( min\{ \mu_{2p}\}, max\{\nu_{2p} \}, max\{w_{2p}
\})$,

here $p =3,4$ $(min\{ \mu_{23},\mu_{24}\}, max\{\nu_{23},\nu_{24}\},
max\{w_{23},w_{24} \})\\ =(min\{ 1,1\}, max\{0.5,0.5\},
max\{0.1,0.1\})=(1 ,0.5 , 0.1)$

Similarly ,we can find  $d_{11}$ and $ d_{31}$ as $d_{11}=(1,
0.7,0.1)$, $d_{31}= (1, 0.8,0.1)$,

$$
[d_{i1}]=\begin{bmatrix} (1,
0.7,0.1)&   \\
(1 ,0.5 , 0.1)& \\
(1, 0.8,0.1) &
\end{bmatrix}
$$

$$
max[s_i]=\begin{bmatrix} 0.95&   \\
0.93& \\
0.92&
\end{bmatrix}
$$
where $s_i=\mu_{11}-\nu_{11}.w_{11}$

\textbf{Step 5:}Finally , we can find an optimum fuzzy set on U as:

$$opt^2_{[d_{i1}]} (U)= \{  u_1/0.95, u_2/0.93 ,  u_3/0.92\}$$

Thus  $u_1$ has the maximum value. Therefore the couple may decide
to buy the car $u_1$.

\section{Conclusion}
In this paper we have redefine the notion of neutrosophic set in a
new way and proposed the concept of neutrosophic soft matrix and
after that different types of matrices in neutrosophic soft theory
have been defiend.then we have introduced some new operations and
properties on these matrices.


\begin{thebibliography}{00}
\bibitem{ans-13} A.Q. Ansaria, R. Biswasb and S.
Aggarwalc, Neutrosophic classifier: An extension of fuzzy classifer,
Applied Soft Computing 13 (2013) 563--573.

\bibitem{ash-02} C. Ashbacher, Introduction to Neutrosophic Logic, American Research Press
Rehoboth 2002.

\bibitem{ata-89} K. Atanassov, G. Gargov, Interval valued intuitionistic fuzzy sets,
Fuzzy Sets Syst. 31 (1989) 343–-349.

\bibitem{ata-99} K.T. Atanassov, Intuitionistic Fuzzy
Sets, Pysica-Verlag A Springer-Verlag Company, New York (1999).

\bibitem{bas-13} T. M.
Basu, N. K. Mahapatra and S. K. Mondal, Intuitionistic Fuzzy Soft
Matrix and Its Application in Decision Making Problems, Annals of
Fuzzy Mathematics and Informatics x/x, (201y), pp.

\bibitem{bor-12} M. J. Borah, T. J. Neog, D. K.
Sut, Fuzzy Soft Matrix Theory And Its Decision Making, International
Journal of Modern Engineering Research, 2 (2012) 121--127.

\bibitem{bro-13} S. Broumi, Generalized Neutrosophic Soft Set
International Journal of Computer Science, Engineering and
Information Technology (IJCSEIT), /2, 2013.

\bibitem{bro-13b}S. Broumi, F. Smarandache, Intuitionistic Neutrosophic Soft
Set, Journal of Information and Computing Science 8/2, (2013),
130--140.

\bibitem{bro-13c} S. Broumi, F. Smarandache and M.
Dhar, On Fuzzy Soft Matrix Based on Reference Function, Information
Engineering and Electronic Business, 2 (2013) 52-59.


\bibitem{bro-14a}S.Broumi,I.Deli,F.Smarandache," Relation on Interval
Neutrosophic Soft Set", Journal of New Results in
Science,2014,submitted.

\bibitem{che-x} Chetia, B. and Das,P.K.,An application of intuitionistic fuzzy
soft matrices in decision making p roblems (communicated)

\bibitem{cag-09ss} N. \c{C}a\u{g}man and S. Engino\u{g}lu, Soft set theory
and uni-int decision making, European Journal of Operational
Research, DOI: 10.1016/j.ejor.2010.05.004.

\bibitem{cag-09sm} N. \c{C}a\u{g}man and S. Engino\u{g}lu,
Soft matrix theory and its decision making, Computers and
Mathematics with Applications 59 (2010) 3308--3314.

\bibitem{cag-14} N. \c{C}a\u{g}man,
Contributions to the theory of soft sets, Journal of New Results in
Science, 4 (2014), 33--41.

\bibitem{cag-12} N. \d{C}a\u{g}man and S. Enginoðlu, Fuzzy soft matrix theory and its
applications in decision making, Iranian Journal of Fuzzy Systems,
9/1 (2012) 109--119.

\bibitem{cag-11e}\d{C}a\u{g}man, N. \c{C}{\i}tak, F. and Engino\u{g}lu, S. FP-soft set theory and its
applications, Annals of Fuzzy Mathematics and Informatics 2/2,
219-226, 2011.

\bibitem{cag-13} N. \d{C}a\u{g}man, S. Karata\c{s},
Intuitionistic fuzzy soft set theory and its decision making,Journal
of Intelligent and Fuzzy Systems 24/4 (2013) 829--836.

\bibitem{deli-12a} \d{C}a\u{g}man, N. Deli, I. Means of FP-Soft Sets and its
Applications, Hacettepe Journal of Mathematics and Statistics, 41/5
(2012) 615--625.

\bibitem{deli-12b}\d{C}a\u{g}man, N. Deli, I. Product of FP-Soft Sets and its
Applications, Hacettepe Journal of Mathematics and Statistics 41/3
(2012) 365 –-374.


\bibitem{del-14} Deli, I."Interval-valued neutrosophic soft sets and its
decision making", http://arxiv.org/abs/1402.3130

\bibitem{del-14c} I .Deli,. S.Broumi ,"neutrosophic soft relation" (communicated)

\bibitem{din-10} B. Dinda and T.K. Samanta, Relations on Intuitionistic Fuzzy Soft
Sets, Gen. Math. Notes, 1/2 (2010) 74--83.

\bibitem{dub-80} Dubois, D. and Prade, H. Fuzzy Set and Systems:
Theory and Applications, Academic Press, New York, 1980.

\bibitem{fen-11} F. Feng, X. Liu , V. L. Fotea, Y. B. Jun, Soft sets and soft rough
sets, Information Sciences 181 (2011) 1125--1137.

\bibitem{fen-10} F. Feng, C. Li,  B. Davvaz, M. Irfan Ali, Soft sets
combined with fuzzy sets and rough sets: a tentative approach, Soft
Computing 14 (2010)  899--911.


\bibitem{gau-93} W.L. Gau, D.J. Buehrer, Vague sets, IEEE Trans. Systems Man
and Cy-bernet, 23 (2) (1993), 610-614.


\bibitem{jia-10} Y. Jiang, Y. Tang, Q. Chen, H. Liu, J.Tang, Interval-valued intuitionistic fuzzy soft sets and their properties,
 Computers and Mathematics with Applications,60
(2010) 906--918.


\bibitem{kal-13} A. Kalaichelvi, P. Kanimozhi, Impact of excessýve television Viewing by children an
analysis using intuitionistic fuzzy soft Matrýces, Int jr. Of
mathematics sciences and applications 3/1 (2013) 103-108.

\bibitem{kha-13} N. Khan, F. H. Khan, G. S. Thakur, Weighted Fuzzy Soft Matrix Theory
and its Decision Making, International Journal of Advances in
Computer Science and Technology   2/10 (2013) 214-218.

\bibitem{kon-09} Z. Kong, L. Gao and L. Wang, Comment on ``A fuzzy soft set theoretic approach
to decision making problems'', J. Comput. Appl. Math. 223 (2009)
540--542.

\bibitem{maj-01b} P.K. Maji,  R.Biswas A.R. Roy, Intuitionistic Fuzzy Soft Sets. The Journal of Fuzzy
Mathematics, 9(3) (2001) 677{-}692.

\bibitem{maj-13} P.K. Maji,  Neutrosophic soft set, Computers and Mathematics with Applications, 45 (2013)
555{-}562.

\bibitem{maj-12}P.K. Maji, A neutrosophic soft set
approach to a decision making problem, Annals of Fuzzy Mathematics
and Informatics, 3/2, (2012), 313--319.

\bibitem{maj-01a} P.K. Maji, R. Biswas and A.R. Roy, Fuzzy soft
sets, Journal of Fuzzy Mathematics, 9(3) (2001) 589{-}602.

\bibitem{mao-13}J. Mao, D. Yao, C. Wang, Group decision making methods based on intuitionistic fuzzy soft
matrices, Applied Mathematical Modelling 37 (2013) 6425-–6436.

\bibitem{mon-13a} J. I. Mondal and T. K. Roy, Some Properties on Intuitionistic
Fuzzy Soft Matrices, International Journal of Mathematics Research
5/2 (2013) 267--276.

\bibitem{mon-14} J. I. Mondal and T. K. Roy, Intuitionistic Fuzzy Soft Matrix Theory and
Multi Criteria in Decision Making Based on T-Norm Operators,
Mathematics and Statistics 2/2 (2014) 55--61.

\bibitem{mon-13b} J. I.
Mondal, T. K. Roy, Theory of Fuzzy Soft Matrix and its Multi
Criteria in Decision Making Based on Three Basic t-Norm Operators,
International Journal of Innovative Research in Science, Engineering
and Technology 2/10 (2013) 5715--5723.

\bibitem{mol-99} D.A. Molodtsov, Soft set theory-first results,
Comput. Math. Appl. 37 (1999) 19{-}31.

\bibitem{muk-08} A. Mukherjee and S.B.Chakraborty, On Intuitionistic fuzzy soft
relations, Bulletin of Kerala Mathematics Association, 5/1 (2008)
35{-}42.

\bibitem{paw-82} Z. Pawlak, Rough sets, Int. J. Comput. Inform. Sci. 11 (1982)
341{-}356.

\bibitem{pei-05} D. Pei and D. Miao, From soft sets to
information systems, In: X. Hu, Q. Liu, A. Skowron, T.Y. Lin, R.R.
Yager, B.Zhang (Eds.), Proceedings of Granular Computing, IEEE 2005,
Volume: 2, pp: 617{-}621.

\bibitem{rab-05} D. Rabounski F. Smarandache L. Borissova Neutrosophic Methods in General Relativity, Hexis, 2005 no:10.

\bibitem{raj-13} P. Rajarajeswari, T. P. Dhanalakshmi, Intuitionistic Fuzzy Soft Matrix Theory And Its Application In
Decision Making, International Journal of Engineering Research
Technology, 2/4 (2013) 1100-1111.

\bibitem{roy-07} A.R. Roy and P.K. Maji, A fuzzy soft set
theoretic approach to decision making problems, J. Comput. Appl.
Math. 203 (2007) 412{-}418.

\bibitem{sam-05} F.Smarandache,"A Unifying Field in Logics. Neutrosophy:
Neutrosophic  Probability, Set and Logic". Rehoboth: American
Research Press,(1998).

\bibitem{sai-13} B.K. Saikia, H. Boruah and P.K. Das, Application of Intuitionistic Fuzzy Soft Matrices in
Decision Making Problems, International Journal of Mathematics
Trends and Technology, 4/11 (2013) 254--265.


\bibitem{sai-14}B.K. Saikia, H. Boruah and P.K. Das, An Appliaction of
Generalized Fuzzy Soft Matrices in Decision Making Problem, IOSR
Journal of Mathematics, 10/1 (2014), PP 33--41.

\bibitem{sez-11} A. Sezgin and A.O. Atag\"{u}n, On operations of
soft sets, Computers and Mathematics with Applications, 61/5 (2011)
1457{-}1467.


\bibitem{som-06} T. Som, On the theory of soft sets, soft relation and fuzzy soft
relation, Proc. of the national conference on Uncertainty: A
Mathematical Approach, UAMA-06, Burdwan, 2006, 1{-}9.

\bibitem{süt-12} D. K. Sut, An application of fuzzy soft relation in decision making
problems, International Journal of Mathematics Trends and Technology
3/2 (2012) 51--54.


\bibitem{wan-05} H. Wang, F. Smarandache, Y.Q. Zhang, R.
Sunderraman, Interval Neutrosophic Sets and Logic: Theory and
Applications in Computing, Hexis; Neutrosophic book series, No: 5,
2005.

\bibitem{yan-09} X. Yang, T.Y. Lin, J. Yang, Y. Li and D. Yu,
Combination of interval-valued fuzzy set and soft set, Comput. Math.
Appl. 58 (2009) 521{-}527.

\bibitem{zad-65}  L.A. Zadeh, Fuzzy Sets, Inform. and Control 8 (1965) 338{-}353.

\bibitem{zha-13} Z. Zhang, C. Wang,
D. Tian, K. Li, A novel approach to interval-valued intuitionistic
fuzzy soft set based decision making, Applied Mathematical
Modelling, xxx (2013) xxx–xxx(ýn press)


\bibitem{zou-08} Y. Zou and Z. Xiao, Data analysis approaches of soft sets under
incomplete information, Knowl. Base. Syst. 21 (2008) 941{-}945.


\end{thebibliography}
\end{document}